\documentclass{amsart}

\usepackage{amsmath, amssymb, mathtools, tikz-cd}

\numberwithin{equation}{section}

\newtheorem{thm}{Theorem}[section]

\newtheorem{cor}[thm]{Corollary}
\newtheorem{prob}[thm]{Problem}

\begin{document}
	
\title{On the AF embedding problem}
\author{Christopher Schafhauser}
\address{Department of Mathematics, 
		University of 	Nebraska -- Lincoln \\ 
		Lincoln, NE 68588, USA.}
\email{cschafhauser2@unl.edu}
\thanks{This research was supported by NSF grant DMS-2400178.}

\begin{abstract}
	This article gives a survey of the AF-embedding problem for $C^*$-algebras.  The problem asks when a given $C^*$-algebra $A$ embeds into an \emph{approximately finite dimensional} (AF) $C^*$-algebra (a direct limit of finite dimensional $C^*$-algebras).  This is largely based on my paper ``Subalgebras of simple AF algebras,'' which shows that separable, exact $C^*$-algebras satisfying the UCT and admitting a faithful, quasidiagonal trace embed into a unital, simple AF algebra.  Some other results are also discussed.
\end{abstract}

\maketitle

\section{Introduction}

A $C^*$-algebra is \emph{approximately finite dimensional} (AF) if it can be written as the closure of an increasing union of finite dimensional $C^*$-algebras.  These $C^*$-algebras were introduced by Bratteli in \cite{Bratteli}, following earlier examples of Glimm (\cite{Glimm}) and Dixmier (\cite{Dixmier}), which in turn were motivated by Murray and von Neumann's construction of the hyperfinite II$_1$ factor $R$ (\cite{Murray-vonNeumann43}).  The structure theory for finite dimensional $C^*$-algebras and morphisms between them led to Bratteli's combinatorial classification of AF algebras in terms of certain directed trees now called Bratteli diagrams, which was later repackaged into Elliott's classification of AF algebras by their $K_0$-group (\cite{Elliott76}).  This classification, combined with the range-of-invariant theorem of Effros, Handelman, and Shen (\cite{EffrosHandelmanShen}), makes AF algebras a very well-understood class of $C^*$-algebras.

Despite the classification of AF algebras, several fundamental questions remain unsolved.  For instance, Blackadar's question on the existence of order two automorphisms of AF algebras with prescribed behavior on $K_0$ is still wide open (\cite[Problem~10.11.3]{Blackadar:kbook}).  This paper is devoted to the problem of finding an abstract characterization of subalgebras of AF algebras. While every von Neumann subalgebra of the hyperfinite II$_1$ factor $R$ is hyperfinite as a consequence of Connes's theorem (\cite{Connes}), it is well-known that AF algebras typically have many non-AF subalgebras.  For a first example, note that if $X$ is a compact metric space, then $X$ is the continuous image of the Cantor space $2^{\mathbb N}$, and hence $C(X)$ embeds into the AF algebra $C(2^{\mathbb N})$, whereas $C(X)$ is itself AF if and only if $X$ is totally disconnected.  A more striking example, and the motivation for the AF embedding problem in general, is given by Pimsner and Voiculescu's embedding of the irrational rotation algebras into AF algebras (\cite{PimsnerVoiculescu:Atheta}), which played a crucial role in the classification of such $C^*$-algebras (\cite{Rieffel:Atheta}).  The following version of the AF embedding problem was popularized by Blackadar and Kirchberg in \cite{BlackadarKirchberg:MathAnn}.

\begin{prob}\label{prob:main}
	If $A$ is a separable, exact, quasidiagonal $C^*$-algebra, is there an AF algebra $B$ and a faithful $^*$-homomorphism $A \rightarrow B$?
\end{prob}

The hypotheses are easily seen to be necessary.  For separability, this is clear since AF algebras are, by definition, separable.  Exactness and quasidiagonality of $A$ together are characterized by the existence of a faithful representation $A \subseteq \mathcal B(\mathcal H)$ together with finite dimensional $C^*$-algebras $F_n$ and completely positive, contractive maps
\begin{equation}
	\begin{tikzcd}
		A \arrow{r}{\phi_n} & F_n \arrow{r}{\psi_n} & \mathcal B(\mathcal H)
	\end{tikzcd}
\end{equation}
such that
\begin{equation}
	\lim_{n \rightarrow \infty} \|\psi_n(\phi_n(a)) - a\| = 0 \quad \text{and} \quad \lim_{n \rightarrow \infty} \|\phi_n(ab) - \phi_n(a) \phi_n(b)\| = 0
\end{equation}
for all $a, b \in A$.  Note that this condition passes to subalgebras as one can simply restrict the domain of the maps $\phi_n$, and it is satisfied for an AF algebra $A$ by taking the $\phi_n$ to be conditional expectations onto finite dimensional subalgebras and the $\psi_n$ to be the inclusions maps.  So AF embeddable $C^*$-algebras are indeed exact and quasidiagonal.  

There have been many partial results verifying the conjecture in large classes of examples; see Section~\ref{sec:history}.  The following two theorems give definitive results under stronger hypotheses.  The first result of Gabe gives a complete solution in the traceless setting, and the second result of myself gives a nearly complete solution in the presence of a faithful trace.

\begin{thm}[{\cite[Corollary~C]{Gabe:GAFA}}]\label{thm:Gabe}
	A separable, exact, traceless $C^*$-algebra embeds into an AF algebra if and only if its primitive ideal space has no non-empty, compact, open sets.
\end{thm}

\begin{thm}[{\cite[Theorem~A]{Schafhauser:AnnMath}}]\label{thm:main}
	A separable, exact $C^*$-algebra satisfying the universal coefficient theorem (UCT) from \cite{RosenbergSchochet:Duke} embeds into a unital, simple AF algebra if and only if it has a faithful, quasidiagonal trace.
\end{thm}

This article is mostly focused on the latter theorem but will discuss some of the techniques and partial results behind the AF embedding problem in general.  Section~\ref{sec:history} surveys several known results and some of the standard techniques related to AF embeddings.  Section~\ref{sec:trace-kernel} outlines the methods behind the proof of Theorem~\ref{thm:main}, particularly the \emph{trace-kernel extension}, and Section~\ref{sec:main} explains how to use these results to obtain Theorem~\ref{thm:main} and some further consequences.  Finally, Section~\ref{sec:further} discusses some results which grew out of the techniques of \cite{Schafhauser:AnnMath}.

\subsection*{Acknowledgments}

Some of the work described here, especially the presentation of the trace-kernel techniques in Section~\ref{sec:trace-kernel}, greatly benefited from my participation in the AIM SQuaRE ``von Neumann algebraic techniques in the classification of $C^*$-algebras'' with Jos\'e Carri\'on, Kristin Courtney, Jamie Gabe, Aaron Tikuisis, and Stuart White during the period 2019--2023.  I am grateful to the American Institute of Mathematics for their support and to my collaborators in this program for several helpful discussions related to the material of this survey over the years.

I acknowledge the support of the Institut Henri Poincar\'e (UAR 839 CNRS-Sorbonne Universit\'e), and LabEx CARMIN (ANR-10-LABX-59-01).
Part of this work was also done during my visit to KU Leuven in November 2026, and I thank the math department for their hospitality.

\subsection*{AI Statement} 

No AI was used in preparing this document.

\section{Previous results}\label{sec:history}

The first non-trivial AF embedding result and the work which led to the abstract AF embedding problem (Problem~\ref{prob:main}) is the classification of the irrational rotation algebras due to Pimsner and Voiculescu in \cite{PimsnerVoiculescu:Atheta} and Rieffel in \cite{Rieffel:Atheta}.  Fix an irrational number $\theta$ and let $A_\theta$ denote the universal $C^*$-algebra generated by unitaries $u$ and $v$ with $u v = e^{2 \pi i \theta} v u$.  In a remarkable tour de force, Pimsner and Voiculescu constructed an embedding $A_\theta \rightarrow B_\theta$, where $B_\theta$ is the unital, simple AF algebra with $K_0(B_\theta) = \mathbb Z + \theta \mathbb Z \subseteq \mathbb R$, with the order inherited from $\mathbb R$ and the unit 1.  This was used by Rieffel to show $K_0(A_\theta) \cong \mathbb Z + \theta \mathbb Z$, which led to classification of the irrational rotation algebras: $A_\theta \cong A_{\theta'}$ if and only if $\theta \equiv \pm \theta'\ (\text{mod}\ \mathbb Z)$.

The basic idea behind the Pimsner--Voiculescu embedding is easy to describe.  While there are no unitaries $u$ and $v$ in a finite dimensional $C^*$-algebra satisfying the relation $uv = e^{2\pi i \theta} vu$, it is possible to satisfy this relation approximately.  In particular, approximate $\theta$ by a rational number $p/q$ where $p, q \in \mathbb Z$, $q \geq 1$, and $p$ and $q$ are relatively prime.  Let $\omega = e^{2 \pi i p/q}$ and define unitaries in $u_0, v_0 \in M_q(\mathbb C)$ by
\begin{equation}
	u_0 = 
	\begin{pmatrix} 
		\omega & & & \\
		& \omega^2 & & \\
		& & \ddots & \\
		& & & \omega^q
	\end{pmatrix}
	\qquad \text{and} \qquad
	v_0 = \begin{pmatrix}
		& & & 1\\
		1 & & & \\
		& \ddots & & \\
		& & 1 &
	\end{pmatrix}.
\end{equation}
Then $u_0 v_0 = e^{2 \pi i p/q} v_0 u_0 \approx e^{2\pi i \theta} v_0 u_0$.

Let $\mathcal Q = \bigotimes_{n=1}^\infty M_n(\mathbb C)$ denote the universal UHF algebra.  Approximating $\theta$ by a sequence of rational numbers $(p_n/q_n)_{n=1}^\infty$ and defining unitaries as above gives sequences of unitaries $(u_n)_{n=1}^\infty, (v_n)_{n=1}^\infty \subseteq \mathcal Q$ which satisfy
\begin{equation}
	\lim_{n \rightarrow \infty} \|u_n v_n - e^{2 \pi i \theta} v_n u_n \| = 0.
\end{equation}
If one could arrange the sequences $(u_n)_{n=1}^\infty$ and $(v_n)_{n=1}^\infty$ to be Cauchy, there would be a unital $^*$-homomorphism $A_\theta \rightarrow \mathcal Q$, which is necessarily injective by the simplicity of $A_\theta$.  However, this is impossible due to the existence of a projection $p \in A_\theta$ of irrational trace (\cite{Rieffel:Atheta}).  Instead, for a suitable sequence of rational approximations $(p_n/q_n)_{n=1}^\infty$ of $\theta$, arising from the continued fraction expansion of $\theta$, Pimsner and Voiculescu construct embeddings
\begin{equation}
	\begin{tikzcd}
		M_{q_n}(\mathbb C) \oplus M_{q_{n+1}}(\mathbb C) \arrow{r} & M_{q_{n+1}}(\mathbb C) \oplus M_{q_{n+2}}(\mathbb C)
	\end{tikzcd}
\end{equation}
such that $(u_n \oplus u_{n+1})_{n=1}^\infty$ and $(v_n \oplus v_{n+1})_{n=1}^\infty$ become Cauchy in the limiting AF algebra $B_\theta$, which yields an embedding $A_\theta \rightarrow B_\theta$.

The example above illustrates the difference between AF embeddability and quasidiagonality.  Writing $\mathcal Q_\omega$ for the norm ultrapower of $\mathcal Q$ with respect to a free ultrafilter $\omega$ on the natural numbers, the sequences $(u_n)_{n=1}^\infty$ and $(v_n)_{n=1}^\infty$ induce unitaries $u, v \in \mathcal Q_\omega$ with $u v = e^{2\pi i \theta} v u$.  This yields a unital embedding $A \rightarrow \mathcal Q_\omega$.  By the nuclearity of $A_\theta$ and the Choi--Effros lifting theorem (\cite{ChoiEffros}), there is a sequence of unital, completely positive maps $\phi_n \colon A_\theta \rightarrow \mathcal Q$ which are isometric and multiplicative in the limit.  Composing these maps with conditional expectations onto sufficiently large matrix subalgebra of $\mathcal Q$ shows $A_\theta$ satisfies Voiculescu's abstract characterization of quasidiagonality from \cite{Voiculescu:QD}.

Another early example of AF embeddability, generalizing the case of the irrational rotation algebras, is due to Pimsner in \cite{Pimsner:AFE}.  Specifically, he showed the following.

\begin{thm}[{\cite[Theorem~9]{Pimsner:AFE}}]\label{thm:Pimsner}
	Let $X$ be a compact metric space and let $\sigma$ be a homeomorphism of $X$.  Then $C(X) \rtimes_\sigma \mathbb Z$ embeds into an AF algebra if and only if there is no open set $U \subseteq X$ with $\overline{\sigma(U)}$ properly contained in $U$.
\end{thm}

As above, it is not too difficult to prove quasidiagonality.  Pimsner shows the condition on open sets implies that for every $x \in X$ and $\epsilon > 0$, there are $x_1, \ldots x_n \in X$ such that $d(x, x_1) < \epsilon$, and $d(\sigma(x_i), x_{i+1}) < \epsilon$ for all $i$, with subscripts computed modulo $n$.  Then there are approximate $^*$-homomorphisms
\begin{equation}
	C(X) \rtimes \mathbb Z \rightarrow M_n(\mathbb C)
\end{equation}
given on $C(X)$ as the direct sum of the point evaluations at $x_1, \ldots, x_n$ and sending the generator of $\mathbb Z$ to the cyclic shift.  Taking a direct sum over finitely many choices of $x$ yields quasidiagonal approximations for $C(X) \rtimes \mathbb Z$.  Again, the difficult part is arranging for these approximations to fit into a Cauchy sequence to obtain a genuine AF embedding instead of an approximate one.

A standard method for verifying AF embeddability has been to produce an AF algebra $B$ and approximately multiplicative, approximately isometric, completely positive maps $\phi_n \colon A \rightarrow B$ which are approximately unitarily equivalent.  If one can show, for example, that there are an increasing sequence of finite subsets $(\mathcal F_n)_{n=1}^\infty$ of $A$ with dense union and a sequence of unitaries $(u_n)_{n=1}^\infty \subseteq B$ such that
\begin{equation}
	\| u_n \phi_n(a) u_n^* - \phi_{n+1}(a) \| < 2^{-n}, \qquad a \in \mathcal F_n,\ n \geq 1,
\end{equation}
then the sequence $(\mathrm{ad}(u_1^*\cdots u_{n-1}^*) \circ \phi_n)_{n=1}^\infty$ is a point-norm Cauchy sequence, and the limiting map $\phi \colon A \rightarrow B$ is an embedding.  This requires a certain ``uniqueness'' result showing that two approximate embeddings are approximately unitarily equivalent.  This is the general strategy behind Dadarlat's AF embedding results in \cite{Dadarlat:MRL, Dadarlat:JFA, Dadarlat:MJM}, for example, where such a uniqueness result is proven up to a matrix amplification via the stable uniqueness theorem (\cite{DadarlatEilers1, DadarlatEilers2, Lin:stable-uniqueness}).

The approach taken to Theorem~\ref{thm:main} is similar in spirit.  Roughly, it uses a uniqueness result for quasidiagonal approximations, showing that two such approximations which approximately agree on $K$-theory and traces are approximately unitarily equivalent, together with an existence result which constructs quasidiagonal approximations which approximate prescribed $K$-theoretic and tracial data.  From here, the quasidiagonal approximations are conjugated into a Cauchy sequence, giving a genuine AF embedding in the limit.

Similar classification results are used in  Gabe's AF embedding theorem in Theorem~\ref{thm:Gabe}.  Here, instead of approximate morphisms into matrices (or AF algebras), Gabe classifies (approximate) morphisms into R{\o}rdam's purely infinite AH algebra $B$ using the classification result from \cite{Gabe:Crelle}.  Writing $B$ as a direct limit of homogeneous algebras $B_n$, $B$ embeds into the (non-separable) AF algebra $\varinjlim B_n^{**}$, and it is easy to construct a separable AF algebra containing $B$ from here.  The classification is in terms of the ideal lattice, and the hypothesis on the primitive ideal space is used to construct a suitable map into the ideal lattice of $B$.

Motivated by Pimsner's Theorem~\ref{thm:Pimsner}, several other AF embedding results for crossed products have also been obtained.  For example, Brown characterized when crossed products of AF algebras by $\mathbb Z$ are AF embeddable in \cite{Brown:JFA} with a condition similar to Pimsner's.  It was also shown that crossed products of UHF algebra by $\mathbb Z^n$ and crossed products of unital, simple A$\mathbb T$ algebra by $\mathbb Z$ are AF embeddable by Brown (\cite{Brown:UHF}) and Matui (\cite{Matui:AT}), respectively.  A simultaneous generalization of these last two results, proving sufficient conditions for AF embeddability for crossed products of AH algebras by $\mathbb Z^n$, was given by Lin in \cite{Lin:AF1, Lin:AF2}.  Some results for crossed products by $\mathbb R$ were obtained by Pimsner in \cite{Pimsner:flows} and Katsura in \cite{Katsura:flows}.

All of the results above used at least some amount of classification techniques.  Through a remarkable application of the classification theory, Ozawa, R{\o}rdam, and Sato showed in \cite{OzawaRordamSato} that all elementary amenable groups are AF embeddable.  This was later extended to all amenable groups by Tikuisis, White, and Winter in \cite{TikuisisWhiteWinter}.  The basic idea is to use classification to show that if $G$ is a countable, discrete, amenable group, then $\big(\big(\bigotimes_G M_{2^\infty} \big) \rtimes G \big) \otimes \mathcal Q$ is an A$\mathbb T$ algebra containing $C^*(G)$ and then to use that that A$\mathbb T$ algebras embed into AF algebras.  A more streamlined proof can be given as a consequence of Theorem~\ref{thm:main}; see Corollary~\ref{cor:group-embedding}.

In another direction, a result essentially due to Spielberg (\cite{Spielberg}) states that if
\begin{equation}
	\begin{tikzcd}
		0 \arrow{r} & I \arrow{r} & A \arrow{r} & B \arrow{r} & 0
	\end{tikzcd}
\end{equation}
is an essential, semisplit extension with vanishing $\mathrm{Ext}$ class such that $I$ and $B$ are AF embeddable, then $A$ is AF embeddable (see \cite[Proposition~2.1]{Dadarlat:cont-fields}).  Spielberg used this to show that all residually stably finite type I $C^*$-algebra are AF embeddable.  It remains open if all stably finite type I $C^*$-algebra are AF embeddable, or even quasidiagonal.  This extension theoretic result was used by Ozawa in \cite{Ozawa:homotopy} to show that, like quasidiagonality, AF embeddability is a homotopy invariant of separable, exact $C^*$-algebras; in particular, the purely infinite $C^*$-algebra $C_0([0, 1)) \otimes \mathcal O_2$ is AF embeddable.  In \cite{Dadarlat:cont-fields}, Dadarlat further used Spielberg's extension result to show AF embeddability for certain continuous fields.

Some further examples are given by various AF embeddability results for graph $C^*$-algebras and related constructions obtained in \cite{Schafhauser:JOT, Schafhauser:BLMS, Schafhauser:JFA, Clark-anHeuf-Sims:JFA, JeongPar:BLMS, AnheufNgSims:JOT}.  Most of these either involve constructing the AF embedding directly at the level of the graph or reduce to the crossed product results of Pimsner and Brown.

At the present time, there is no clear approach to the general AF embedding problem, but there are two very general results worth mentioning.  Blackadar and Kirchberg showed in \cite{BlackadarKirchberg:MathAnn} that separable, nuclear, quasidiagonal $C^*$-algebras admit a generalized inductive limit structure, given by a sequence of finite dimensional $C^*$-algebras and completely positive, contractive maps which, in a suitable sense, are multiplicative in the limit.  Also, Dadarlat proved the following general ``external'' AF approximation property in \cite[Theorem~6]{Dadarlat:external}.  If $A \subseteq \mathcal B(H)$ is a separable, exact, quasidiagonal $C^*$-algebra with $A \cap \mathcal K(\mathcal H) = 0$ and $\mathcal H$ separable, then there is a (not necessarily increasing) sequence of finite dimensional $C^*$-algebras $F_n \subseteq \mathcal B(\mathcal H)$ such that $A \subseteq \overline{\bigcup_{n=1}^\infty F_n}$; of course, if one can choose the sequence $(F_n)_{n=1}^\infty$ to be increasing, this would give an affirmative solution to Problem~\ref{prob:main}.

\section{The trace-kernel extension}\label{sec:trace-kernel}

Returning to the setting of Theorem~\ref{thm:main}, as suggested in the previous section, the basic strategy is to classify approximate embeddings $A \rightarrow B$ for a suitable AF algebra $B$.  Instead of approximate morphisms, it is technically easier to work with morphisms into ultrapowers.  To this end, for a unital $C^*$-algebra $B$, let $B_\omega$ denote the norm ultrapower of $B$.  Note that for a $C^*$-algebra $A$, any sequence of approximate embeddings $A \rightarrow B$ induces an embedding $A \rightarrow B_\omega$.

The following result from \cite{Schafhauser:AnnMath} is the formal statement of the classification of approximate AF embeddings mentioned in the previous section.  When $A$ is nuclear, an analogous result for much more general codomains $B$ is given in \cite[Theorem~9.1]{CarrionGabeEtal:ActaMath}.

\begin{thm}[{\cite[Propositions~4.2 and~4.3]{Schafhauser:AnnMath}}]\label{thm:approx-embedding}
	Let $A$ be a unital, separable, exact $C^*$-algebra satisfying the UCT with a faithful, quasidiagonal trace $\tau$, and let $B$ be a unital, simple, $\mathcal Q$-stable AF algebra.  Then every unital, trace-preserving group homomorphism $K_0(A) \rightarrow K_0(B_\omega)$ is induced by a unital, nuclear, trace-preserving embedding $A \rightarrow B_\omega$, which is unique up to unitary equivalence.
\end{thm}

The quasidiagonality condition on the trace $\tau$ and the exactness of $A$ give the existence of a nuclear, trace-preserving $^*$-homomorphism $\psi_0 \colon A \rightarrow \mathcal Q_\omega$ such that $\tau_{\mathcal Q_\omega} \circ \psi_0 = \tau$, where $\tau_{\mathcal Q_\omega}$ denotes the unique trace on $\mathcal Q_\omega$.  Since $\tau$ is assumed to be faithful, $\psi_0$ is necessarily an embedding.  Further, as every non-zero corner of $\mathcal Q_\omega$ is isomorphic to $\mathcal Q_\omega$, we may assume $\psi_0$ is unital.  As $B$ is unital and $\mathcal Q$-stable, there is a unital embedding $\mathcal Q \rightarrow B$, and composing $\psi$ with the ultrapower of this embedding yields a unital, nuclear embedding $\psi \colon A \rightarrow B_\omega$ such that $\tau_{B_\omega} \circ \psi = \tau$, where $\tau_{B_\omega}$ is the unique trace on $B_\omega$.  This is the starting point in the proof of Theorem~\ref{thm:approx-embedding}.

The first key observation, motivated by Matui and Sato's work on the Toms--Winter conjecture in \cite{MatuiSato:Acta,MatuiSato:Duke} and the quasidiagonality theorem of Tikuisis, White, and Winter in \cite{TikuisisWhiteWinter} (see also \cite{Gabe:QD, Schafhauser:Crelle}), is that $B_\omega$ admits a von Neumann algebra quotient.  More precisely, let $\pi$ denote the GNS representation of the unique trace on $B$ and note that $\pi(B)''$ is isomorphic to the hyperfinite II$_1$ factor $R$.  A standard consequence of Kaplansky's density theorem shows that $\pi$ induces a surjective $^*$-homomorphism $q \colon B_\omega \rightarrow R^\omega$, giving an exact sequence
\begin{equation}\label{eq:trace-kernel}
	\begin{tikzcd}
		0 \arrow{r} & J \arrow{r}{j} & B_\omega \arrow{r}{q} & R^\omega \arrow{r} & 0,
	\end{tikzcd}
\end{equation}
called the \emph{trace-kernel extension} of $B$.  The algebra $J = \ker(q) \unlhd B_\omega$ is known as the \emph{trace-kernel ideal}.

A standard consequence of Connes's theorem (\cite{Connes}) states that unital, nuclear $^*$-homomorphisms $A \rightarrow R^\omega$ are determined up to unitary equivalence by their trace.  Further, since the unitary group of $R^\omega$ is path-connected, any unitary in $R^\omega$ lifts to a unitary in $B_\omega$.  Hence any unital, nuclear, trace-preserving $^*$-homomorphism $A \rightarrow B_\omega$ is unitarily equivalent to such a $^*$-homomorphism $\phi$ with $\mathrm{im}(\phi - \psi) \subseteq J$.  So Theorem~\ref{thm:approx-embedding} amounts to classifying $J$-perturbations of $\psi$ by their induced map on $K_0$.  Taking some care with the non-separability of $J$, such maps are classified up to stabilization and homotopy by their difference $[\phi, \psi]$ in the Kasparov--Skandalis group $KK_{\mathrm{nuc}}(A, J)$, defined in \cite{Skandalis}, refining Kasparov's $KK$-theory from \cite{Kasparov}.  The basic goal becomes to compute $KK_{\mathrm{nuc}}(A, J)$ in terms of $K_0$ and to replace stabilized homotopy with unitary equivalence.

For the first part, this is an immediate consequence of the UCT.  First one verifies $K_1(B_\omega) = 0$ and $K_1(R^\omega) =0 $.  Then the six-term exact sequence in $K$-theory, applied to the trace-kernel extension in \eqref{eq:trace-kernel}, gives rise to an exact sequence
\begin{equation}
	\begin{tikzcd}
		0 \arrow{r} &  K_0(J) \arrow{r}{j_0} \arrow{r} & K_0(B_\omega) \arrow{r}{q_0} & K_0(R^\omega) \arrow{r}{\partial_0} & K_1(J) \arrow{r} & 0.
	\end{tikzcd}
\end{equation}
By the classification of projections in the II$_1$ factor $R^\omega$, the trace $\tau_{R^\omega}$ induces an isomorphism $\hat\tau_{R^\omega} \colon K_0(R^\omega) \rightarrow \mathbb R$.  It follows that $K_0(J)$ is the kernel of the map $\hat\tau_{B_\omega} \colon K_0(B_\omega) \rightarrow \mathbb R$ and $K_1(J) = 0$.  Also, using the $\mathcal Q$-stability of $B$, one can show $K_0(J) \cong K_0(J) \otimes \mathbb Q$ and hence $K_*(J)$ is divisible.  So the UCT yields an isomorphism $KK_{\mathrm{nuc}}(A, J) \cong \mathrm{Hom}(K_0(A), K_0(J))$.  Repackaging this, there is a bijection from the set of $\kappa \in KK_{\mathrm{nuc}}(A, J)$ with $\kappa_0([1_A]) = 0$ and the set of unital, trace-preserving group homomorphisms $K_0(A) \rightarrow K_0(B_\omega)$ defined by sending $\kappa$ to $\kappa_0 + K_0(\psi)$.  So Theorem~\ref{thm:approx-embedding} is a consequence of the following result.

\begin{thm}\label{thm:KK-classification}
	Let $A$ be a unital, separable, exact $C^*$-algebra with a faithful, quasidiagonal trace $\tau$, let $B$ be a unital, simple, $\mathcal Q$-stable AF algebra,
	and let $\psi \colon A \rightarrow B_\omega$ be a unital, nuclear, trace-preserving $^*$-homomorphism.  If $\kappa \in KK_{\mathrm{nuc}}(A, J)$ with $\kappa_0([1_A]) = 0$, then there is a unital, nuclear, trace-preserving $^*$-homomorphism $\phi \colon A \rightarrow B_\omega$ such that $\mathrm{im}(\phi - \psi) \subseteq J$ and $[\phi, \psi] = \kappa$; further, this $\phi$ is unique up to unitary equivalence.
\end{thm}

The result above was proven in a somewhat ad hoc way in \cite{Schafhauser:AnnMath}, reusing material from \cite{Schafhauser:Crelle}; the discussion given here grew out of \cite{CarrionGabeEtal:ActaMath}.  The underlying mechanism is the notion of an \emph{absorbing representation} (\cite{Thomsen:rep}), which can be viewed as an abstraction of the conclusion of Voiculescu's non-commutative Weyl--von Neumann theorem (\cite{Voiculescu:rep, Arveson}).  For a separable, stable $C^*$-algebra $I$, let $\mathcal M(I)$ denote the multiplier algebra of $I$.  A unital $^*$-homomorphism $\bar\psi \colon A \rightarrow \mathcal M(I)$ is called \emph{(unitally, nuclearly) absorbing} if for every unital, nuclear $^*$-homomorphism $\rho \colon A \rightarrow \mathcal M(I)$, there is a sequence of $1 \times 2$ unitary matrices $(u_n)_{n=1}^\infty$ over $\mathcal M(I)$ such that
\begin{equation}
	u_n (\rho(a) \oplus \bar\psi(a)) u_n^* - \bar\psi(a) \in I
\end{equation}
and such that this sequence converges to 0 in norm; usually one requires this for all \emph{weakly nuclear} $\rho$, but using the exactness of $A$, any weakly nuclear map $\rho$ is nuclear by \cite[Proposition~3.2]{Gabe:JNCG}.

Using the Hjelmborg--R{\o}rdam characterization of stability in \cite{HjelmborgRordam} and strict comparison for $B_\omega$, one can show that $J$ is \emph{separably stable} in the sense of \cite{Schafhauser:AnnMath}: every separable subalgebra of $J$ is contained in a separable, stable subalgebra of $J$.  This is shown in \cite[Propositions~3.2 and~3.3]{Schafhauser:AnnMath} (see also \cite[Lemma~6.11]{CarrionGabeEtal:ActaMath}).  To access absorption, one should replace the trace-kernel ideal with a sufficiently large separable subalgebra at this stage, but for the sake of exposition, we will ignore this step and treat $J$ as if it was separable.  Since $J$ is an ideal of $B_\omega$, there is a natural unital $^*$-homomorphism $\lambda \colon B_\omega \rightarrow \mathcal M(J)$.  Using Elliott and Kucerovsky's characterization of absorption (\cite[Theorem~6]{ElliottKucerovsky}) together with the later work of Kucerovsky and Ng (\cite{KucerovskyNg}) and Ortega, Perera, and R{\o}rdam (\cite{OrtegaPereraRordam}) on the \emph{corona factorization property}, one can show that the map $\bar\psi = \lambda \circ \psi \colon A \rightarrow \mathcal M(J)$ is absorbing.  The key point here is the $\mathcal Q$-stability of $J$, which follows from the $\mathcal Q$-stability of $B$.

The existence portion of Theorem~\ref{thm:KK-classification} follows from a general \emph{$KK$-existence theorem} (see \cite[Proposition~2.6]{Schafhauser:AnnMath}, for example).  Roughly, using that $\bar\psi$ is unital, nuclear, and absorbing, one finds a unital, nuclear  $\bar\phi \colon A \rightarrow \mathcal M(J)$, which is necessarily absorbing, such that $\mathrm{im}(\bar \phi - \bar \psi) \in J$ and $[\bar\phi, \bar\psi] = \kappa$.  Then 
\begin{equation}
	\bar\phi(A) \subseteq \bar\psi(A) + J \subseteq \lambda(B_\omega),
\end{equation} 
and since $\lambda$ is injective, $\bar\phi$ factors as $\bar \lambda = \lambda \circ \phi$ for the desired map $\phi$.

The uniqueness portion of Theorem~\ref{thm:KK-classification} is more subtle.  Given $\phi_1, \phi_2 \colon A \rightarrow \mathcal M(B)$  as in the statement, $[\phi_1, \psi] = [\phi_2, \psi]$ implies $[\phi_1, \phi_2] = 0$.  Also, $\bar\phi_i = \lambda \circ \phi_i$, $i = 1, 2$, is absorbing since $\lambda \circ \psi$ is.  Then $\bar\phi_1$ and $\bar \phi_2$ absorb each other and hence are approximately unitarily equivalent via a sequence of unitaries $\mathcal M(J)$.  The key is to use that fact that $[\bar \phi_1, \bar \phi_2] = 0$ to show that these unitaries can be taken in the minimal unitization $J$, and hence in $\lambda(B_\omega)$; using that $\lambda$ is isometric, $\phi_1$ and $\phi_2$ are approximately unitarily equivalent and so a reindexing argument shows that they are unitarily equivalent.  This was done in \cite[Proposition~2.7]{Schafhauser:AnnMath} by using the Dadarlat--Eilers stable uniqueness theorem (\cite[Theorem~3.10]{DadarlatEilers1}) to prove the result up to a matrix amplification and then using $\mathcal Q$-stability to hide this matrix amplification.  A more conceptual proof is now possible through Szab\'o's recent \emph{$KK$-uniqueness theorem} in \cite[Corollary~4.8]{Szabo}.

\section{Theorem~\ref{thm:main} and consequences}\label{sec:main}

The bulk of the work behind the proof of Theorem~\ref{thm:main} is contained in Theorem~\ref{thm:approx-embedding}, which gives a classification of approximate embeddings $A \rightarrow B$ for $A$ as in the statement of Theorem~\ref{thm:main} and suitable AF algebras $B$.  Via standard intertwining arguments (see \cite[Section~5]{Schafhauser:AnnMath}, for example), one obtains the following classification result for genuine embeddings from the classification of approximate embeddings.  As suggested in Section~\ref{sec:history}, the basic strategy is to conjugate the approximate embeddings into a Cauchy sequence and take the limit.

\begin{cor}[{\cite[Theorem~D]{Schafhauser:AnnMath}}]\label{cor:embedding}
	If $A$ is a unital, separable, exact $C^*$-algebra satisfying the UCT with a faithful, quasidiagonal trace $\tau$ and $B$ is a unital, simple, $\mathcal Q$-stable AF algebra, then every unital, trace-preserving group homomorphism $K_0(A) \rightarrow K_0(B)$ is induced by a unital, nuclear, trace-preserving embedding $A \rightarrow B$, which is unique up to approximate unitary equivalence.
\end{cor}

To obtain an AF embedding as in Theorem~\ref{thm:main}, all that remains is to find a unital, simple, $\mathcal Q$-stable AF algebra such that there is a unital, trace-preserving group homomorphism $K_0(A) \rightarrow K_0(B)$ since Corollary~\ref{cor:embedding} then provides the AF embedding $A \rightarrow B$.  This is in fact easy from known results: particularly, Elliott's predecessor of the Effros--Handelman--Shen theorem (\cite{EffrosHandelmanShen}) in \cite{Elliott:total-order}, which shows that all countable, totally ordered abelian groups are dimension groups.  One takes $B$ to be the AF algebra with $K_0(B)$ isomorphic to the $\mathbb Q$-linear span of $\hat\tau(K_0(A)) \subseteq \mathbb R$, so that $\hat\tau$ gives a unital, trace-preserving map $K_0(A) \rightarrow K_0(B)$.

Some further applications of Corollary~\ref{cor:embedding} are given in \cite[Section~6]{Schafhauser:AnnMath}.  I highlight two of them here.  The following gives a sharper version of the AF embeddability of amenable group $C^*$-algebras obtained in \cite{OzawaRordamSato, TikuisisWhiteWinter}.  In addition to Corollary~\ref{cor:embedding}, the main tools are \cite{HigsonKasparov, Tu}, showing that $C^*(G)$ satisfies the UCT, and \cite{Luck}, showing that range of the trace on $K_0$ is contained in $\mathbb Q$.

\begin{cor}[{\cite[Theorem~B]{Schafhauser:AnnMath}}]\label{cor:group-embedding}
	If $G$ is a countable, discrete, amenable group, then $C^*(G)$ admits a unital, trace-preserving embedding into $\mathcal Q$.
\end{cor}

The following gives a general result for AF embeddability of crossed products.  This is immediate from Theorem~\ref{thm:main}, again using \cite{HigsonKasparov, Tu} to verify the UCT condition.

\begin{cor}[{\cite[Theorem~C]{Schafhauser:AnnMath}}]\label{cor:crossed-product}
	If $X$ is a compact metric space, $G$ is a countable, discrete, amenable group, and $G \curvearrowright X$ is an action such that $X$ admits a faithful, $G$-invariant probability measure, then $C(X) \rtimes G$ embeds into an AF algebra.
\end{cor}

Without the existence of a faithful trace, the AF embedding problem for crossed products $C(X) \rtimes G$ is still open.  Pimsner's Theorem~\ref{thm:Pimsner} gives a complete solution for $G = \mathbb Z$, but the AF embedding problem is unsolved for $C(X) \rtimes \mathbb Z^2$, for example, even when $X$ is the Cantor set.  This old problem was recently restated as \cite[Problem~XI]{SchafhauserTikuisisWhite}.  In a different direction, Rainone and I gave a complete characterization of quasidiagonality for $A \rtimes \mathbb Z$ when $A$ in an AH algebra of real rank zero (\cite[Theorem~1.1]{RainoneSchafhauser:AdvMath}), but it is not known if this also characterizes AF embeddability for such crossed products.

\section{Further work}\label{sec:further}

The trace-kernel methods outlined in Section~\ref{sec:trace-kernel} have proven to be a powerful tool for problems related to the classification and structure theory for $C^*$-algebras.  We mention a few here.

First, a significant refinement of Theorem~\ref{thm:approx-embedding} (and Corollary~\ref{cor:embedding}) was given in \cite{CarrionGabeEtal:ActaMath}, which classifies morphisms between unital, separable, simple, nuclear, $\mathcal Z$-stable $C^*$-algebras satisfying the UCT by a suitable combination of $K$-theoretic and tracial invariants.  This lead to a new proof of the classification theorem for such $C^*$-algebras.  We refer to the introduction of \cite{CarrionGabeEtal:ActaMath} for a more thorough discussion.  The methods were also extended in \cite{Schafhauser24} as a step towards a classification theorem which does not assume the UCT.  Remarkably, Bouwen and Gabe showed in \cite{BouwenGabe} that the trace-kernel techniques can be adapted to the non-tracial setting by replacing the trace with a suitable state.  This was used to obtain a new proof of the Kirchberg--Phillips theorem in the spirit of the proof of the stably finite classification theorem given in \cite{CarrionGabeEtal:ActaMath}.

Aside from classification, some structural results for non-nuclear $C^*$-algebras have been obtained through trace-kernel methods.  For example, Kunnawalkam Elayavalli and I proved results on the exponential rank and exponential length of certain $C^*$-algebras in \cite{Tarski}.  This was used to show that the reduced free group $C^*$-algebras $C^*_r(F_n)$ are mutually not elementarily equivalent.  The classification of AF embeddings in Corollary~\ref{cor:embedding} was also used to obtain new examples of MF groups in \cite{Schafhauser:GGD}.  Finally, \cite{K-stability} uses trace-kernel methods to give new examples of $K$-stable $C^*$-algebras.

\end{document}